\documentclass[11pt,a4paper]{article}
\usepackage{amssymb}

\newtheorem{deff}{Definition}[section]
\newtheorem{lemma}[deff]{Lemma}
\newtheorem{theorem}[deff]{Theorem}
\newtheorem{corollary}[deff]{Corollary}

\newtheorem{fact}[deff]{Fact}

\newtheorem{em-example}[deff]{Example}
\newtheorem{em-def}[deff]{Definition}        
\newtheorem{em-remark}[deff]{Remark}         
\newtheorem{em-question}[deff]{Question}

\newenvironment{definition}{\begin{em-def} \em  }{ \end{em-def}}
\newenvironment{remark}{\begin{em-remark} \em }{\end{em-remark}}
\newenvironment{question}{\begin{em-question}\em }{\end{em-question}}
\newenvironment{proof}{\noindent {\it Proof}.}{\hfill\mbox{$\square$}\smallskip}

\newcommand{\EOP}{\hfill\mbox{$\square$}\smallskip}
\def\ker{\mathop{\rm ker}}

\def\T{{\mathbb T}}
\def\Z{{\mathbb Z}}
\def\N{{\mathbb N}}
\def\R{{\mathbb R}}
\def\Q{{\mathbb Q}}
\def\P{{\mathbb P}}
\def\cont{\mathfrak c}

\begin{document}
\title{{Quasi-convex density and determining subgroups of compact abelian groups}}
\author{Dikran Dikranjan\footnote{Dipartimento di Matematica e Informatica, Universit\`{a} di Udine, Via delle Scienze  206, 33100 Udine, Italy; {\em e-mail\/}: 
{\tt dikranja@dimi.uniud.it}; the first author was partially supported by MEC.MTM2
006-02036 and FEDER FUNDS.
}\mbox{\hskip40pt } and \mbox{\hskip40pt }Dmitri Shakhmatov\footnote{%
Graduate School of Science and
 Engineering, Division of Mathematics, Physics and Earth Sciences, Ehime University, Matsuyama 790-8577,  Japan; {\em e-mail\/}: {\tt dmitri@dpc.ehime-u.ac.jp};
 the second author was partially supported by the Grant-in-Aid for Scientific 
Research no.~19540092 by the Japan Society for the Promotion of Science (JSPS).}
}
\date{\small{\em Dedicated to W. Wistar Comfort on the occasion of his 75th anniversary}}
\maketitle

\begin{abstract} 
For an abelian topological group $G$, let $\widehat{G}$ denote the
dual group of all continuous characters endowed with the compact open
topology.  Given a closed subset $X$ of an infinite compact abelian
group $G$ such that $w(X)<w(G)$, and an open neighbourhood $U$ of $0$
in $\T$, we show that $|\{\chi\in\widehat{G}: \chi(X)\subseteq
U\}|=|\widehat{G}|$.  (Here, $w(G)$ denotes the weight of $G$.)  A
subgroup $D$ of $G$ determines $G$ if the map $r:
\widehat{G}\to\widehat{D}$ defined by $r(\chi)=\chi\restriction_D$ for
$\chi\in \widehat{G}$, is an isomorphism between $\widehat{G}$ and
$\widehat{D}$.  We prove that
$$
w(G)=\min\{|D|: D\mbox{ is a subgroup of $G$ that determines }G\}
$$ 
for every infinite compact abelian group $G$. In particular, an
infinite compact abelian group determined by a countable subgroup is
metrizable.  This gives a negative answer to questions of Comfort,
Hern\' andez, Macario, Raczkowski and Trigos-Arrieta from \cite{CRT,
CRT1, HMT}.  As an application, we furnish a short elementary proof of
the result from \cite{HMT} that a compact abelian group $G$ is
metrizable provided that every dense subgroup of $G$ determines $G$.
\end{abstract}

\bigskip
All topological groups are assumed to be Hausdorff, and all topological spaces are assumed to be Tychonoff. As usual, $\R$ denotes the group of real numbers (with the usual topology), $\Z$ denotes the group of integer numbers, $\T=\R/\Z$ denotes the circle group (with the usual topology), $\N$ denotes the set of natural numbers, $\P$ denotes the set of prime numbers, $\omega$ denotes the first infinite cardinal, and $w(X)$ denotes the weight of a space $X$.  If $A$ is a subset of a space $X$, then $\overline{A}$ denotes the closure of $A$ in $X$.

\section{Preliminaries and background}

In this section we give necessary definitions and collect five facts that will be needed later. These facts are either known or part of the folklore. However, to make this manuscript self-contained, we provide their proofs in Section \ref{facts} for the reader's convenience.

For spaces $X$ and $Y$,  we denote by $C(X,Y)$ the space of all continuous functions from $X$ to $Y$ endowed with the 
{\em compact open topology\/}, that is, the topology generated by the family 
$$
\{[K,U]: K \mbox{ is a compact subset of } X \mbox{ and } U \mbox{ is an open subset of } Y\}
$$ 
as a subbase, where
$$
[K,U]=\{g\in C(X,Y): g(K)\subseteq U\}.
$$

\begin{fact}\label{weight:C(C,T)}
If $X$ is a compact space and $Y$ is a space, then $w(C(X,Y))\le w(X)+w(Y)+\omega$.
\end{fact}

For a topological group $G$, we denote by $\widehat{G}$
 the Pontryagin-van Kampen dual of $G$, namely 
 the group of all continuous characters $\chi:G\to\T$ endowed 
with the compact open topology.  Clearly, $\widehat{G}$ is a closed subgroup of $C(G,\T)$. In particular, a base of neighborhoods of $0$ in $\widehat{G}$ is given by the sets 
\begin{equation}
\label{definition:of:W:K:U}
W(K,U)= \{\chi\in\widehat{G}: \chi(K)\subseteq U\}= [K,U]\cap \widehat{G},
\end{equation}
where $K$ is a compact subset of $G$ and $U$ is an open neighbourhood of $0$ in $\T$. 

We identify $\T=\R/\Z$ with the real interval $(-1/2, 1/2]$ in the obvious way, and write
$$
\T_+=\{x\in \T: -1/4\leq x \leq 1/4 \}.
$$

\begin{definition}
Let $G$ be an abelian topological group. 
\begin{itemize}
\item[(i)]
For $E\subseteq G$ and $A \subseteq \widehat{G}$, define the
{\em polars\/}
$$
E^\triangleright=\{\chi\in \widehat{G}\; |\; \chi(E)\subseteq \T_+\}
\;\mbox{ and }\;
A^\triangleleft=\{x\in {G}\; |\; \chi(x)\in \T_+ \mbox{ for all }
\chi\in A\}.
$$ 
\item[(ii)]
A set $E\subseteq G$ is said to be {\em quasi-convex\/} if $E=E^{\triangleright\triangleleft}$.
\item[(iii)]
The {\em quasi-convex hull\/} $Q_{G} (E)$ of $E\subseteq G$ is the smallest quasi-convex set of $G$ containing $E$. 
\item[(iv)]
Following \cite{DDL,dL}, we will say that 
 $E\subseteq G$ is {\em qc-dense\/} (an abbreviation for {\em quasi-convexly dense\/})
provided that $Q_G(E)=G$, or equivalently, if $E^{\triangleright}=\{0\}$.
\end{itemize}
\end{definition}

Obviously, $E\subseteq E^{\triangleright\triangleleft}$. Therefore, a set $E\subseteq G$
is quasi-convex if and only if for every $x\in G\setminus E$ there exists $\chi\in E^\triangleright$ such that $\chi(x)\not\in \T_+$. 

The notion of quasi-convexity was introduced by Vilenkin \cite{Vil} as a natural counterpart for topological groups of 
the fundamental notion of convexity from  the theory of topological vector spaces (we refer the reader to \cite{B,Diss} for additional information).

\begin{fact}
\label{finite:sum:is:qc-dense}
Suppose that $U$ is an open neighbourhood of $0$ in $\T$ and $X$ is a compact subset of a topological group $G$ such that $W(X,U)=\{0\}$. Then:
\begin{itemize}
\item[\rm{(i)}]
 there exists $n\in\N$ such that the sum 
$$
K_n=(X\cup\{0\})+(X\cup\{0\})+\dots+(X\cup\{0\})
$$
 of $n$ many copies of the set $X\cup\{0\}$ is qc-dense in $G$;
\item[\rm{(ii)}] the  subgroup of $G$ generated by $X$ must be dense in $G$.
\end{itemize}
\end{fact}

Item (i) of our next fact can be found in \cite{dL,DDL}.

\begin{fact}\label{images} Let $f:G\to H$ be a continuous homomorphism of topological abelian groups. Then:
\begin{itemize}
\item[\rm{(i)}] $f(Q_G(X))\subseteq Q_H(f(X))$ for every subset $X$ of $G$.
\item[\rm{(ii)}]
If $f(G)$ be dense in $H$ and $X$ is a qc-dense subset of $G$, then $f(X)$ is qc-dense in $H$. 
\end{itemize}
\end{fact}

\begin{definition}
Following \cite{CRT,CRT1}, we say that a subgroup $D$ of an abelian group $G$
{\em determines\/} $G$ if the restriction homomorphism $r: \widehat{G}\to \widehat{D}$ (defined by $r(\chi)=\chi\restriction_D$ for $\chi\in \widehat{G}$)
is an isomorphism between the topological groups $\widehat{G}$ and $\widehat{D}$.
\end{definition}

 This notion is relevant 
to
extending the Pontryagin-van Kampen duality to non-locally compact groups
\cite{CM, Diss}. 
Indeed, if $G$ is locally compact and abelian, then every subgroup $D$ that determines $G$ must be dense in $G$; in particular, no proper locally compact subgroup of $G$ can determine $G$. 
 (We note that the original definition in \cite{CRT,CRT1} assumed upfront that $D$ is dense in $G$.) When $D$ is dense in $G$, the restriction homomorphism $r: \widehat{G}\to \widehat{D}$  is always a continuous isomorphism.

The ultimate connection between the notions of determined subgroup and qc-density is established in the next fact. This fact is a particular case of a more general fact stated without proof
(and in equivalent terms) in \cite[Remark 1.2(a)]{CRT1} and \cite[Corollary 2.2]{HMT}.
\begin{fact}  
\label{connection:between:determination:and:qc-density}
For a subgroup $D$ of a compact abelian group $G$  the following conditions are equivalent:
\begin{itemize}
\item[\rm{(i)}]
$D$ determines $G$;
\item[\rm{(ii)}]
there exists a compact subset of $D$ which is qc-dense in $G$.
\end{itemize}
\end{fact}

\begin{definition}
According to \cite{CRT, CRT1}, a topological group $G$ is said to be {\em  determined\/} if every dense subgroup of $G$ determines $G$.
\end{definition}

Chasco 
\cite[Theorem 2]{CM}
and Au\ss enhofer \cite[Theorem 4.3]{Diss} proved that all metrizable abelian groups are determined.  Comfort,
Raczkowski and  Trigos-Arrieta established the following amazing inverse of this theorem for compact groups: Under the Continuum
Hypothesis CH, every determined compact abelian group is metrizable (\cite[Corollary 4.9]{CRT} and  \cite[Corollary 4.17]{CRT1}). Quite recently, Hern\'{a}ndez,  Macario and Trigos-Arrieta removed the assumption of CH from their result \cite[Corollary 5.11]{HMT}. We note that this theorem becomes  an immediate  consequence of our main result, see Corollary  \ref{determined:groups:corollary}.

\begin{fact}
\label{quotient:of:determined}
{\rm \cite[Corollary 3.15]{CRT1}}
If $f:G\to H$ is a continuous surjective  homomorphism of compact abelian groups and $G$ is determined, then $H$ is determined as well.
\end{fact}

The following question remained the last principal unsolved problem in the theory of compact determined groups:

\begin{question}\label{Ques1}{\rm  \cite[Question 7.1]{CRT},  \cite[Question 7.1]{CRT1}, \cite[Question 5.12]{HMT}} 
\begin{itemize}
\item[(a)] Is there a compact group $G$ with a countable dense subgroup $D$ such that $w(G) > \omega$ and $D$ determines $G$? 
\item[(b)] What if $G=\T^\kappa$? 
\end{itemize}
\end{question}

We completely resolve this question in Corollary \ref{countable:subgroups:determining:G}.
In fact, we even solve the most general version of this question with $\omega$ replaced by an arbitrary cardinal, see Corollary \ref{size:of:determined:subgroup}.

\section{Main results}

\begin{definition}
\label{restriction:map}
If $X$ is a subset of a compact abelian group $G$, then $r^G_X: \widehat{G}\to C(X,\T)$ denotes the ``restriction map'' defined by $r^G_X(\chi)=\chi\restriction_X$ for 
$\chi\in\widehat{G}$. 
\end{definition}

Observe that $C(X,\T)$ is a topological group and  $r^G_X$ is a continuous group homomorphism. 

We refer the reader to (\ref{definition:of:W:K:U}) for the definition of $W(X,U)$.
\begin{theorem}
\label{scattered:subsets:are:big}
Let $X$ be a  closed subset of  an infinite compact abelian group $G$ such that  $w(X)<w(G)$. Then for every open neighbourhood $U$ of $0$ in $\T$ one has $|W(X,U)|=|\widehat{G}|$. 
\end{theorem}

\begin{proof}
Consider first the case when $w(X)<\omega$. Then $X$ must be finite.
Note that the  set $W(X,U)$ is an open neighborhood of $0$ in the initial topology ${\mathcal T}$ of $\widehat{G}$ with respect to the family $\{\eta_x: x\in X\}$
of evaluation characters $\eta_x: \widehat{G} \to \T$ defined by $\eta_x(\pi)=\pi(x)$ for every $\pi \in \widehat{G}$. Since topologies generated by characters are totally bounded, finitely many translations of $W(X,U)$ cover the whole group $\widehat{G}$. Since $\widehat{G}$ is infinite, this yields   $|W(X,U)|=|\widehat{G}|$. 

From now on we will assume that $w(X)\ge \omega$.  The inequality $|W(X,U)|\le |\widehat{G}|$ being trivial, it suffices to check that $|\widehat{G}|\le |W(X,U)|$.

Let $r^G_X$ be the map from Definition \ref{restriction:map}, and let  $H=r^G_X(\widehat{G})$. Note that  $\ker r^G_X \subseteq W(X,U)$, so  $|\ker r^G_X|\le |W(X,U)|$. If $|\ker r^G_X|=|\widehat{G}|$, we  are done.  Assume now that $|\ker r^G_X|<|\widehat{G}|$.  Since $\widehat{G}$ is infinite, we obtain
\begin{equation}
\label{0:eq:1}
|\widehat{G}|=|\widehat{G}/\ker r^G_X|= |r^G_X(\widehat{G})|=|H|.
\end{equation}

Let $N$ be the subgroup of $H$ generated by the open subset
$[X,U]\cap H$ of $H$. Then $N$ is a clopen subgroup of $H$, so the index of $N$ in $H$ cannot exceed $w(H)$, which gives 
\begin{equation}
\label{0:eq:2}
  |H|=|N|+|H/N|\le |N|+w(H)\le |[X,U]\cap H|+\omega +w(H).
\end{equation}
Since  $w(H)\le w(C(X,\T))\le w(X)$ by Fact \ref{weight:C(C,T)}, and $w(X)+\omega=w(X)$ by our assumption, we obtain from (\ref{0:eq:2}) that
\begin{equation}
\label{0:eq:3}
|H|\le |[X,U]\cap H|+w(X).
\end{equation}
As $w(X)<w(G)=|\widehat{G}|=|H|$ by (\ref{0:eq:1}),  and $|H|=|\widehat{G}|\ge \omega$, from (\ref{0:eq:3}) it follows that
\begin{equation}
\label{0:eq:4}
|[X,U]\cap H|=|H|.
\end{equation}
Finally, note that $[X,U]\cap H = r^G_X(W(X,U))$, which yields that
\begin{equation}
\label{0:eq:5}
|[X,U]\cap H|=|r^G_X(W(X,U))|\le |W(X,U)|.
\end{equation}
Combining (\ref{0:eq:1}), (\ref{0:eq:4}) and (\ref{0:eq:5}), we obtain the inequality $|\widehat{G}|\le |W(X,U)|$.
\end{proof}

\begin{corollary}
\label{scattered:subspaces:are:big}
If a  closed subspace $X$ of  an infinite compact abelian group $G$ is qc-dense in $G$, then  $w(X)=w(G)$.
\end{corollary}
\begin{proof}
Let $U$ be an open neighbourhood of $0$ in $\T$ such that $U\subseteq \T_+$.
Since $X$ is qc-dense in $G$, we have $W(X,U)\subseteq X^{\triangleright}=\{0\}$. Now Theorem \ref{scattered:subsets:are:big} yields  $w(X)\ge w(G)$. The reverse inequality $w(X)\le w(G)$ is trivial.
\end{proof}

Our next corollary constitutes a major breakthrough in the theory of compact determined groups. 

\begin{corollary}
\label{size:of:determined:subgroup}
If a subgroup $D$ of an infinite compact abelian group $G$ determines $G$, then  $|D|\ge w(G)$.
\end{corollary}
\begin{proof}
According to Fact \ref{connection:between:determination:and:qc-density}, $D$ contains a compact subset $X$ that is qc-dense in $G$, so
$|D|\ge |X|\ge w(X)$ (see, for example, \cite[Theorem 3.1.21]{Eng}). Finally, $w(X)=w(G)$ by Corollary \ref{scattered:subspaces:are:big}.
\end{proof}

Even the particular case of Corollary \ref{size:of:determined:subgroup} provides a complete answer to Question \ref{Ques1}:

\begin{corollary}
\label{countable:subgroups:determining:G}
A compact abelian group determined by a countable  subgroup is metrizable.
\end{corollary}

\begin{corollary}
\label{determined:groups:corollary} {\rm \cite[Corollary 5.11]{HMT}}
Every determined compact abelian  group is metrizable.
\end{corollary} 
\begin{proof}
Assume that $G$ is a non-metrizable determined compact abelian group. Then $w(G)\ge\omega_1$, and so we can find a continuous surjective group homomorphism  $h: G\to K=H^{\omega_1}$, where $H$ is either $\T$ or $\Z(p)$ for some prime number $p$
(see, for example, \cite[Theorem 5.15 and Discussion 4.14]{CRT1}).
As a continuous homomorphic image of the determined group $G$, the group $K$ is determined by Fact \ref{quotient:of:determined}.  
Since $K$ is separable (see, for example, \cite[Theorem 2.3.15]{Eng}), 
there exists a countable dense subgroup $D$ of $K$. Since $K$ is determined, we conclude that $D$ must determine $K$. Therefore, $K$ must be metrizable by Corollary \ref{countable:subgroups:determining:G}, a contradiction.
\end{proof}

Useful properties of determined groups can be found in \cite{CT}.

A {\em super-sequence\/} is a non-empty compact Hausdorff space $X$ with at most one non-isolated point $x^*$ \cite{DS}. When $X$ is infinite, we will call $x^*$ the {\em limit\/} of $X$ and say that $X$ {\em converges to $x^*$\/}.  Observe that a convergent sequence is a countably infinite super-sequence.

Au\ss{}enhofer \cite{Diss} essentially proved that every infinite compact metric abelian group has a qc-dense sequence converging to 0.\footnote{This is an immediate consequence of \cite[Theorem 4.3 or Corollary 4.4]{Diss}. In fact, a more general statement immediately follows from these results: Every dense subgroup $D$ of a compact metric abelian group $G$ contains a sequence converging to 0 which is qc-dense in $G$.}  Our next theorem extends this result to all compact groups by replacing converging sequences with super-sequences. 

\begin{theorem}
\label{theorem:2}
Every infinite compact abelian group contains a  qc-dense super-sequence converging to $0$.
\end{theorem}

The proof of Theorem \ref{theorem:2} is postponed until Section \ref{Section:3}.

\begin{corollary}
\label{getting:a:small:determined:subgroup}
Every infinite compact abelian group $G$ has a (dense) subgroup $D$ which determines $G$ such that $|D|\le w(G)$. 
\end{corollary}

\begin{proof}
Apply Theorem \ref{theorem:2} to find a super-sequence $X$ that is qc-dense in $G$. Let $D$ be the 
subgroup of $G$ 
generated by $X$.  Clearly, 
$|X|=w(X)\le w(G)$. Since $G$ is infinite, $w(G)$ must be infinite,  and therefore $|D|\le \omega+ |X|\le w(G)$. Finally, $D$ determines $G$ by Fact \ref{connection:between:determination:and:qc-density}.
\end{proof}

Our next corollary provides another major advance in the theory of compact determined groups:

\begin{corollary}
\label{dd:G}
If $G$ is an infinite compact abelian group, then
$$
w(G)=\min\{|D|: D \mbox{ is a subgroup of } G \mbox{ that determines } G\}.
$$
\end{corollary}

\begin{proof}
Combine Corollaries \ref{size:of:determined:subgroup} and \ref{getting:a:small:determined:subgroup}.
\end{proof}

We have been kindly informed by Chasco that our next corollary was independently proved by Bruguera and Tkachenko:

\begin{corollary}
Every infinite compact abelian group $G$ contains a proper (dense) subgroup $D$ which determines $G$. 
\end{corollary}

\begin{proof}
Let $D$ be a subgroup of $G$ as in the conclusion of Corollary \ref{getting:a:small:determined:subgroup}. Since $G$ is an infinite compact group, we have 
$|D|\le w(G)<2^{w(G)}=|G|$. Therefore, $D$ must be a proper subgroup of $G$.
\end{proof}

A subspace $X$ of a topological group $G$ {\em topologically generates $G$\/} if $G$ is the smallest  closed subgroup of $G$ that contains $X$. 
\begin{remark}
\label{remark:1}
\begin{itemize}
\item[\rm{(i)}] Item (ii) of Fact \ref{finite:sum:is:qc-dense} can be restated as follows: {\em A qc-dense subset of a compact abelian group $G$ topologically generates $G$.\/} 
Therefore, for a subset $X$ of a compact abelian group $G$, one has the following implications:
\begin{equation}
\label{some:implications}
X \mbox{ is dense in } G 
\longrightarrow
X \mbox{ is qc-dense in } G 
\longrightarrow
X \mbox{ topologically generates } G. 
\end{equation} 
\item[\rm{(ii)}] 
The first arrow in (\ref{some:implications}) cannot be reversed. Indeed, take any qc-dense sequence $S$ in $\T$ 
(see Lemma \ref{qc-dense:in:T} for an example of such a sequence). Clearly, $S$ is not dense in $\T$.

\item[(iii)] 
The last arrow in (\ref{some:implications}) cannot be reversed either.  Indeed, it follows from the results in \cite{DS} that $\T^\cont$ contains a converging sequence (i.e., countably infinite super-sequence) topologically generating $\T^\cont$. This sequence, however, cannot be qc-dense in $\T^\cont$  by Corollary \ref{scattered:subspaces:are:big}.
\end{itemize}
\end{remark}

\label{remark:HM}
According to a well-known result of Hofmann and Morris 
\cite{HM, HMbook},
every compact group $G$ contains a super-sequence topologically generating $G$.  (See also \cite{Sh} for a ``purely topological'' proof of this result based on Michael's selection theorem.) The emphasized text in Remark \ref{remark:1}(i) allows us to  conclude that Theorem \ref{theorem:2} implies the particular case of the theorem of Hofmann and Morris for compact {\em abelian\/} groups $G$.  As  it  was demonstrated in Remark \ref{remark:1}(iii),  a (super-)sequence topologically generating a compact abelian group $G$ need not be qc-dense in $G$. Therefore, the conclusion of our Theorem \ref{theorem:2} is {\em formally stronger\/} than that of (the abelian case of) the result of Hofmann and Morris. 

\section{Characterization of qc-dense subsets and determining subgroups of compact abelian groups in 
terms of $C(X,\T)$}

We refer the reader to (\ref{definition:of:W:K:U}) and Definition \ref{restriction:map} for notations used in our next theorem.
\begin{theorem}
\label{near:characterization}
For a closed subset $X$ of a compact abelian group $G$ the following conditions are equivalent:
\begin{itemize}
\item[\rm{(i)}] $W(X,U)=\{0\}$ for some open neighbourhood $U$ of 0 in $\T$;
\item[\rm{(ii)}] $r^{G}_X$ is an isomorphism between the topological groups $\widehat{G}$ and $H=r^G_X(\widehat{G})$.
\end{itemize}
\end{theorem} 
\begin{proof}
(i)$\to$(ii)
Let $U$ be as in (i).  Since $\ker r^G_X \subseteq W(X,U)=\{0\}$, we conclude that $r^G_X$ is an injection. Since $X$ is compact, 
$$
\{r^G_X(0)\}=r^G_X(\{0\})=r^G_X(W(X,U))=H\cap \{g\in C(X,\T): g(X)\subseteq U\}
$$
is an open subset of $H$. Since $H$ is a subgroup of $C(X,\T)$, we conclude that $H$ is discrete. Therefore, $r^G_X$ is an open map onto its image.

(ii)$\to$(i) The assumption from item (ii) implies that $H$ is a discrete subgroup of $C(X,\T)$. Hence, we can find $n\in\N$, compact subsets $K_0,\dots,K_n$ of $X$ and open neighbourhoods $U_0,\dots, U_n$ of $0$ in $\T$ such that
\begin{equation}
\label{1:eq:1}
H\cap \bigcap_{i\le n} \{g\in C(X,\T): g(K_i)\subseteq U_i\}=\{r^G_X(0)\}.
\end{equation}
Define $U=\bigcap_{i\le n} U_i$. Now equation (\ref{1:eq:1}) yields $W(X,U)=\{0\}$. 
\end{proof}

\begin{corollary}
\label{qc-denseness:implies:isomorphism}
If a closed subset $X$ of a compact abelian group $G$ is qc-dense in $G$, then $r^G_X$ is  an isomorphism between the topological groups  $\widehat{G}$ and $r^G_X(\widehat{G})$.
\end{corollary}
\begin{proof}
Choose an open neighbourhood $U$ of $0$ with $U\subseteq \T_+$.
Since $X$ is qc-dense in $G$, we have $W(X,U)\subseteq X^\triangleright=\{0\}$,
and we can apply Theorem \ref{near:characterization} to this $U$.
\end{proof}

\begin{corollary}
\label{isomorphism:implies:qc-density:of:some:combination}
Let $X$ be a closed subset  of a compact abelian group $G$ such that
$r^G_X$ is an isomorphism between the topological groups $\widehat{G}$ and $r^G_X(\widehat{G})$.  Then there exists $n\in \N$ such that  the sum 
$$
K_n=(X\cup\{0\})+(X\cup\{0\})+\dots+(X\cup\{0\})
$$
 of $n$ many copies  of the set $X\cup\{0\}$ is (compact and) qc-dense in $G$.
\end{corollary}
\begin{proof}
Apply Theorem \ref{near:characterization} to find an open neighbourhood  $U$ of 0 in $\T$ as in item (i) of this theorem. Then apply Fact \ref{finite:sum:is:qc-dense}(i) to obtain the required $n\in\N$.
\end{proof}

\begin{corollary}
For a subgroup $D$ of a compact abelian group $G$ the following conditions are equivalent:
\begin{itemize}
\item[\rm{(i)}] $D$ determines $G$;
\item[\rm{(ii)}] there exists a compact set $X\subseteq D$ such that $r^{G}_X$ is an isomorphism between the topological groups $\widehat{G}$ and $r^G_X(\widehat{G})$.
\end{itemize}
\end{corollary}
\begin{proof}
(i)$\to$(ii) Since $D$ determines $G$, there exists a compact set $X\subseteq D$ which is qc-dense in $G$ (Fact \ref{connection:between:determination:and:qc-density}). Then $r^{G}_X$ is  an isomorphism between the topological groups $\widehat{G}$ and $r^G_X(\widehat{G})$ (Corollary \ref{qc-denseness:implies:isomorphism}).

(ii)$\to$(i) Let $X$ be as in item (ii). Apply Corollary \ref{isomorphism:implies:qc-density:of:some:combination} to get $n\in\N$ and $K_n$ as in the conclusion of this corollary. Clearly, $K_n$ is a compact subset of $D$. Since $K_n$ is qc-dense in $G$, $D$ determines $G$ by Fact \ref{connection:between:determination:and:qc-density}.
\end{proof}

\section{Proof of Theorem \ref{theorem:2}}
\label{Section:3}

We start with a partial inverse of Fact \ref{images}(ii). 

\begin{lemma}\label{3space}
Suppose that $f:G\to H$ is a continuous surjective homomorphism of compact abelian groups  and $X$ is a subset of $G$ such that $X\cap \ker f$ is qc-dense in  $\ker f$. Then $X$ is qc-dense in $G$ if and only if $f(X)$ is qc-dense in $H$.  
\end{lemma}

\begin{proof} If $X$ is  qc-dense in $G$, then $f(X)$ is qc-dense in $H$ by Fact \ref{images}(ii). 

Assume that $f(X)$ is qc-dense in $H$. Let $\chi\in X^\triangleright$. Since
$$
Q_G(X)\supseteq Q_G(X\cap \ker f)\supseteq Q_{\ker f}(X\cap \ker f)=\ker f,
$$
one has $\chi\in (\ker f)^\triangleright$. Since $\ker f$ is a subgroup and $\T_+$ contains no non-trivial subgroups, this yields that $\chi$ vanishes on $\ker f$. Thus, $\chi$ factorizes as $\chi= \xi\circ f$, where $\xi \in \widehat H$. Now $\chi \in X^\triangleright$ obviously yields $\xi \in f(X)^\triangleright$.  
As $f(X)$ is qc-dense in $H$ by the hypothesis,  this yields $\xi \in f(X)^\triangleright=\{0\}$. Hence $\xi =0$, and so $\chi=0$ as well. Therefore, $X^\triangleright=\{0\}$, and thus $X$ is qc-dense in $G$.
\end{proof}

The following definition is an adaptation to the abelian case of \cite[Definition 4.5]{DS}:
\begin{definition}
\label{definition:of:fan}
Let $\{G_i:i\in I\}$ be a family of abelian topological groups. For every $i\in I$, let $X_i$ be a subset of $G_i$. Identifying each $G_i$ with a subgroup of the direct product 
$G=\prod_{i\in I} G_i$ in the obvious way, define $X=\bigcup_{i\in I} X_i\cup\{0\}$, where  $0$ is the zero element of $G$. We will call $X$ the {\em fan\/} of the family $\{X_i:i\in I\}$ and will denote it by $\mathrm{fan}_{i\in I}(X_i, G_i)$. 
\end{definition}

The proof of the following lemma is straightforward.
\begin{lemma}
\label{fan:of:sequences:is:a:super-sequence}
Let $\{G_i:i\in I\}$ be a family of abelian topological groups. For every $i\in I$, let $X_i$ be a sequence converging to $0$ in $G_i$. Then  $\mathrm{fan}_{i\in I}(X_i, G_i)$ is a super-sequence in $G=\prod_{i\in I} G_i$ converging to $0$. 
\end{lemma}

\begin{lemma}\label{products:have:qc-dense}
Let $\{G_i:i\in I\}$ be a family of abelian topological groups. For each $i\in I$ let $X_i$ be a qc-dense subset of $G_i$. Then $X=\mathrm{fan}_{i\in I}(X_i, G_i)$ is qc-dense in $G=\prod_{i\in I} G_i$. 
\end{lemma}
\begin{proof} Let $\chi:G\to\T$ be a non-trivial continuous character. There exist a non-empty finite
subset $J$ of $I$ and a family $\{\chi_j\in\widehat{G_j}:j\in J\}$ such that $\chi(g)=\sum_{j\in J} \chi_j(g(j))$ for $g\in G$. Since $J\not=\emptyset$, we can fix $j_0\in J$. Since $X_{j_0}$ is qc-dense in $G_{j_0}$, there exists $x\in X_{j_0}\subseteq X$
such that $\chi_{j_0}(x)\not\in\T_+$. Finally, note that 
$$
\chi(x)=\sum_{j\in J} \chi_j(x(j))
=
\chi_{j_0}(x(j_0))+\sum_{j\in J\setminus\{j_0\}} \chi_j(x(j))
=
\chi_{j_0}(x)+\sum_{j\in J\setminus\{j_0\}} \chi_j(0)
=
\chi_{j_0}(x)\not\in\T_+.
$$
Therefore, $\chi\not\in X^\triangleright$. This gives $X^\triangleright=\{0\}$, and so $X$ is qc-dense in $G$.
\end{proof}

Our next three lemmas are particular cases of a general result of  Au\ss{}enhofer quoted in the text preceding Theorem \ref{theorem:2}.
Her proof relies on Arzela-Ascoli theorem and inductive construction, so the qc-dense sequence she constructs in her proof is  ``generic''. To keep this manuscript self-contained, we provide   ``constructive'' examples of ``concrete'' qc-dense sequences in the circle group (Lemma \ref{qc-dense:in:T}), 
in the group of $p$-adic integers (Lemma \ref{Z_p:has:qc-dense}) and in the dual group of the rationals equipped with the discrete topology (Lemma \ref{sequence:in:Qhat}). 

\begin{lemma}\label{qc-dense:in:T}
Let $T=\left\{\frac{1}{2n}:n\in\N\right\}\cup\{0\}\subseteq\R$, and let $\varphi: \R\to \R/\Z$ be the natural quotient map. Then $\varphi(T)$ is  a converging to $0$ sequence in 
$\T=\R/\Z$
  that is qc-dense in $\T$.
\end{lemma}

\begin{proof}
Let $\chi\in\widehat{\T}$ be a non-zero character.  Then there exists $m\in\Z\setminus\{0\}$ such that $\chi(x)=mx$ for all $x\in\T$. 
Let $n=|m|$. Then  $\frac{1}{2n}\in T$ and so $x=\varphi\left(\frac{1}{2n}\right)\in\varphi(T)$.
Moreover, $\chi(x)=mx=\varphi\left(\frac{m}{2n}\right)\not\in\T_+$, which shows that $\chi\not\in
\varphi(T)^{\triangleright}$. Hence $\varphi(T)^{\triangleright}=\{0\}$, and so $\varphi(T)$ is qc-dense in $\T$.
\end{proof} 

\begin{lemma}\label{Z_p:has:qc-dense}
For every prime number $p$, the group $\Z_p$ of $p$-adic integers contains a qc-dense sequence converging to $0$.
\end{lemma}
\begin{proof}
Recall that the family $\{p^n\Z_p:n\in\N\}$ consisting of clopen subgroups of $\Z_p$ forms a basis of neighbourhoods of $0$. Therefore,
$S = \{kp^n: n\in \N, 1\le k\le p-1\}\cup\{0\}$ is a sequence converging to $0$ in $\Z_p$.
Let us show that $S$ is qc-dense in $\Z_p$. To this end take a non-zero character $\chi: \Z_p\to \T$.
Since $\Z_p$ is a zero-dimensional compact group, its image $\chi(\Z_p)$ under the continuous homomorphism $\chi$ must be a closed zero-dimensional subgroup of $\T$.
In particular, $\chi(\Z_p)\ne \T$. Being a proper closed subgroup of $\T$, 
$\chi(\Z_p)$ must be finite. It follows that $\ker \chi$ is a clopen subgroup of
$\Z_p$, and so there exists $n\in \N$ such that  $\ker \chi =p^n\Z_p$.
Hence $\chi(\Z_p)\cong \Z_p/p^n\Z_p\cong \Z(p^n)$. Therefore, $\chi(1)=\frac{m}{p^n}$ for some $m$ coprime to $p$. Choose $k$ with $1\le k\le p-1$ and such that:
\begin{itemize}
\item[(a)] $km\equiv \frac{p-1}{2} (mod \; p)$, if $p>2$;
\item[(b)] $k=1$ if $p=2$.
\end{itemize}
Then  $\chi(kp^{n-1})=\frac{p-1}{2p} \not\in \T_+$, 
in case (a). Otherwise, $\chi(2^{n-1})=\frac{1}{2} \not\in \T_+$, in case (b). 
In both cases, $\chi(kp^{n-1})\not\in \T_+$, so $\chi\not\in S^\triangleright$. This proves that $S^\triangleright=\{0\}$. Therefore, 
$S$ is qc-dense in $\Z_p$.
\end{proof}

\begin{lemma}
\label{sequence:in:Qhat}
Let $\Q$ be the group of rational numbers with the discrete topology. Then $\widehat{\Q}$ contains a qc-dense sequence converging to $0$.
\end{lemma}
\begin{proof} 
By Lemma \ref{Z_p:has:qc-dense},  for every prime number $p$ the group $\Z_p$ contains a qc-dense sequence $S_p$ 
converging to $0$. By Lemma \ref{products:have:qc-dense}, $S=\mathrm{fan}_{p\in \P}(S_p, \Z_p)$ is qc-dense  in $K=\prod_{p\in \P} \Z_p$. In view of Lemma \ref{fan:of:sequences:is:a:super-sequence}, $S$ is a sequence converging to $0$ in $K$.

Define   $v=\{1_p\}_{p\in\P}\in K$, where each $1_p$ is the identity of $\Z_p$. Let $N=\R \times K$  and $u=(1,v)\in N$.
Then the cyclic subgroup  $\langle u\rangle$  of  $N$  is discrete and  $G=N/\langle u\rangle$ is isomorphic to $\widehat{\Q}$  \cite[\S 2.1]{DM}, so we will identify  $\widehat{\Q}$ with the quotient $G=N/\langle u\rangle$.  Define $H=N/(\Z\times K)$ and note that $H\cong \T$. Let  $\theta: N \to G= N/\langle u\rangle$  and  $\psi: N\to H$ 
be  the (continuous) quotient homomorphisms. Since $\ker \theta=\langle u\rangle\subseteq \Z\times K=\ker \psi$, there exists a (unique) continuous  surjective group homomorphism $f: G\to H$  such that $\psi=f\circ \theta$.
In particular, $\ker f = \theta(\ker \psi)=\theta(\Z\times K)$. Since $\Z\times K=\langle u\rangle +( \{0\}\times K)$, we have 
\begin{equation}
\label{equation:eight}
\ker f= \theta(\Z\times K)=
\theta(\langle u\rangle +(\{0\}\times K))=
\theta(\langle u\rangle)+\theta(\{0\}\times K)=
\theta(\{0\}\times K).
\end{equation}

Let $T$ and $\varphi$ be as in Lemma \ref{qc-dense:in:T}.  Since $T$ is a sequence converging to $0$ in $\R$, and $S$ is a sequence converging to $0$ in $K$,  it follows that 
$X=\theta(T \times \{0\})\cup \theta(\{0\}\times S)$ is a sequence converging to $0$ in $G$. 

Note that $\varphi(T)=\psi(T\times\{0\})=f(\theta(T\times\{0\}))$ is qc-dense in $f(G)=H\cong \T$ by 
Lemma \ref{qc-dense:in:T}. Since $\theta(T\times\{0\})\subseteq X$, we conclude that  $f(X)$ is qc-dense in $H$ as well. Since $S$ is qc-dense in $K$, 
$\theta(\{0\}\times S)$ is qc-dense in $\theta(\{0\}\times K)=\ker f$ by Fact \ref{images}(ii) and (\ref{equation:eight}).
From $\theta(\{0\}\times S)\subseteq X$, we conclude that $X\cap \ker f$ is qc-dense in $\ker f$. It now follows from Lemma \ref{3space} that $X$ is qc-dense in $G$.
\end{proof}

The next lemma is probably known, but we include its proof for the reader's convenience. 

\begin{lemma}
\label{quotients:of:torsion:free} Every infinite compact abelian group of weight $\kappa$ is isomorphic to a 
quotient group of the group 
$\widehat{\Q}^\kappa \times \prod_{p\in \P} \Z_p^{\kappa}$. 
\end{lemma}

\begin{proof}  Let $H$ be an infinite compact abelian group 
such that $w(H)=\kappa$. Clearly, $\kappa$ is infinite and $X=\widehat H$ is a discrete abelian group of size $\kappa$ \cite[Theorem (24.15)]{hr}. Let $Y = X \oplus \bigoplus_{\kappa}  (\Q \oplus \Q/\Z)$. 
By \cite[Theorem 24.2]{fuchs} there exists a divisible abelian group $D$ containing $Y$ such that no proper subgroup of $D$ containing $Y$ is divisible. 
According to the text immediately following \cite[Theorems 24.2]{fuchs},  $r_0(D)=r_0(Y)$ and $r_p(D)=r_p(Y)$ for every prime $p$, where 
$r_0(N)$ and $r_p(N)$ denote the free-rank and the $p$-rank of an abelian group $N$, 
respectively (see, for example, \cite[\S 16]{fuchs}). Since $r_0(Y)=\kappa$ and $r_p(Y)=\kappa$ for every prime $p$, we conclude that 
\begin{equation}
\label{equation:NEW}
D\cong  \bigoplus_{\kappa} (\Q \oplus \Q/\Z)\cong \left(\bigoplus_{\kappa} \Q\right) \oplus \bigoplus_{p\in \P}
 \left(\bigoplus_{\kappa}\Z(p^\infty)\right)
\end{equation} 
by the structure theorem for divisible abelian groups (see \cite[Theorem 23.1]{fuchs}).
Consider the compact group $G= \widehat D$. By (\ref{equation:NEW}), $G \cong \widehat{\Q}^\kappa \times \prod_{p\in \P} \Z_p^{\kappa}$. 
According to \cite[Theorem (24.5)]{hr}, $H\cong \widehat{X}\cong G/X^\perp$, where $X^\perp=\{\chi\in \widehat{D}: \chi(X)=\{0\}\}$. Therefore, $H$ is isomorphic to a quotient 
group of $\widehat{\Q}^\kappa \times \prod_{p\in \P} \Z_p^{\kappa} \cong G $. 
 \end{proof}

\bigskip
\noindent
\textbf{Proof of Theorem \ref{theorem:2}}: 
Let $H$ be an infinite compact abelian group. Define $\kappa=w(H)$ and $G= \widehat{\Q}^\kappa \times \prod_{p\in \P} \Z_p^{\kappa}$. By Lemma \ref{quotients:of:torsion:free}
there exists a surjective continuous homomorphism $f: G \to H$. Clearly, $G=\prod_{i\in I} G_i$, where $|I|=\kappa$ and each $G_i$ is either 
$\widehat{\Q}$ or $\Z_p$ for a suitable $p\in \P$. By Lemmas \ref{sequence:in:Qhat} and \ref{Z_p:has:qc-dense}, for every $i\in I$ there exists a sequence 
$X_i$ converging to $0$ which is qc-dense in $G_i$. Applying Lemmas \ref{fan:of:sequences:is:a:super-sequence} and \ref{products:have:qc-dense}, 
we conclude that $X=\mathrm{fan}_{i\in I}(X_i, G_i)$ is a super-sequence in $G$ converging to $0$ such that $X$ is qc-dense in $G$. Since $f:G\to H$ is a surjection, $S=f(X)$ is qc-dense in $H$ by Fact \ref{images}(ii). Being the image of a super-sequence $X$ in $G$ converging to $0$,  $S$ is super-sequence in $H$ converging to $0$ \cite[Fact 4.3]{DS}.  Finally, $|S|\le|X|\le\omega\cdot |I|=\kappa$.
\EOP

\section{Proofs of Facts
\ref{weight:C(C,T)},   \ref{finite:sum:is:qc-dense}, \ref{images}, \ref{connection:between:determination:and:qc-density} and \ref{quotient:of:determined}
}\label{facts}

\noindent
\textbf{Proof of Fact \ref{weight:C(C,T)}}:
Indeed, let $\mathcal{B}$ and $\mathcal{C}$ be bases for $X$ and $Y$ respectively such that
$|\mathcal{B}|\le w(X)$ and $|\mathcal{C}|\le w(Y)$. It suffices to show that the family
$\mathcal{E}=\{[\overline{V},U]:V\in\mathcal{B}, U\in\mathcal{C}\}$ is a subbase for the topology of $C(X,Y)$, since this would imply $w(C(X,Y))\le |\mathcal{B}\times\mathcal{C}|\le w(X)+w(Y)+\omega$. 
Indeed, assume that $f\in C(X,Y)$, $K$ is a compact subset of $X$, $O$ is an open subset of $Y$ and $f\in [K,O]$. Let $x\in K$. Since $f(x)\in O$ and $\mathcal{C}$ is a base for $Y$, 
there exists $U_x \in \mathcal{U}$ such that $f(x)\in U_x\subseteq O$.
Using continuity of $f$ and regularity of $X$ one can find $V_x\in\mathcal{B}$ such that $f(\overline{V_x})\subseteq U_x$.  
Since $K$ is compact, $K\subseteq \bigcup_{x\in F} V_x$ for some finite 
subset $F$ of $K$. Now $f\in\bigcap_{x\in F} [\overline{V_x},U_x]\subseteq [K,O]$.
\EOP

\bigskip\noindent
\textbf{Proof of Fact \ref{finite:sum:is:qc-dense}}:
(i)
There exists $n\in \N$ such that 
\begin{equation}
\label{eq:Vn}
V_n=\{x\in \T: kx\in \T_+\mbox{ for all } k=1,2,\ldots, n\}\subseteq U.
\end{equation}

Let $\chi\in K_n^\triangleright$. Fix $x\in X$. Let $k=1,2,\ldots, n$ be arbitrary. Since $0\in X\cup\{0\}$, one has  $kx\in K_n$, and so  $k\chi(x)=\chi(kx)\in \T_+$.
This yields $\chi(x)\in V_n\subseteq U$ by (\ref{eq:Vn}). Since $x\in X$ was chosen arbitrarily, it follows that $\chi\in W(X,U)$. Since $W(X,U)=\{0\}$, this gives $\chi=0$. Therefore, $K_n^\triangleright=\{0\}$,  and so $K_n$ is qc-dense in $G$.

(ii) Suppose that the smallest subgroup $N$ of $G$ containing $X$ is not dense in $G$. Then we can choose $\chi\in\widehat{G}$
such that $\chi(\overline{N})=\{0\}$ and $\chi(y)\not=0$ for some $y\in\widehat{G}\setminus \overline{N}$. So $\chi\in W(X,U)$ and yet $\chi\not=0$, in contradiction with our assumption.
\EOP

\bigskip\noindent
\textbf{Proof of Fact \ref{images}}:
(i) Assume that $x\in Q_G(X)$ but $f(x)\not \in Q_H(f(X))$. Then there exists $\xi\in \widehat{H}$
such that $\xi(f(X))\subseteq \T_+$ and $\xi(f(x))\not \in \T_+$. Then $\chi=\xi \circ f\in \widehat{G}$
and $\chi(X)\subseteq \T_+$, while $\chi(x)\not\in \T_+$. Therefore, $x\not \in    Q_G(X)$, a 
contradiction. 

(ii) By our assumption,  $Q_G(X)=G$. Therefore, $f(G)=f(Q_G(X))\subseteq Q_H(f(X))$ by item (i). Since $ Q_H(f(X))$ must be closed in $H$ and $f(G)$ is dense in $H$, this yields $Q_H(f(X))=H$,  that is,  $f(X)$ is qc-dense in $H$. 
\EOP

\bigskip\noindent
\textbf{Proof of Fact \ref{connection:between:determination:and:qc-density}}:
Since $\widehat{G}$ is discrete, (i) is equivalent to $\widehat{D}$ being discrete.
Since $\widehat{D}$ carries the compact open topology,  this is equivalent to  having  $W(K,U)=\{0\}$ for an
appropriate pair of a compact subset $K$ of $D$ and an open neighborhood $U$ of $0$ in $\T$.  Having this in mind, we are going to prove that (i) and (ii) are equivalent.

(ii)$\to$(i) 
Suppose that $K$ is a compact subset of $D$ that is qc-dense in $G$. Take an
open neighbourhood of $0$ in $\T$ with $U\subseteq \T_+$. Then $W(K,U)\subseteq K^\triangleright =\{0\}$, which gives $W(K,U)=\{0\}$. Thus, (i) holds.

(i)$\to$(ii) By our assumption, there exist a compact subset $X$ of $D$ and an open neighborhood $U$ of $0$ in $\T$ such that $W(X,U)=\{0\}$. Let $K_n$ be as in Fact \ref{finite:sum:is:qc-dense}. Then $K_n\subseteq D$ is  compact and qc-dense in $G$.
\EOP

\bigskip\noindent
\textbf{Proof of Fact \ref{quotient:of:determined}}:
Let $D$ be a dense subgroup of $H$. Then $f^{-1}(D)$ is a dense subgroup of $G$.
Since $G$ is determined, by Fact \ref{connection:between:determination:and:qc-density} we can find a compact subset  $X$
of $f^{-1}(D)$ that is qc-dense in $G$. 
Then $f(X)$ 
is a compact subset of $D$ which is qc-dense in $H$ by 
Fact \ref{images}(ii). Applying Fact \ref{connection:between:determination:and:qc-density} once again, we conclude that $D$ determines $H$.
\EOP

\bigskip
\noindent
{\bf Acknowledgement.\/} 
The authors would like to thank Lydia Au\ss{}enhofer for attracting their attention to Theorem 4.3 and Corollary 4.4 of her manuscript
\cite{Diss}, and G\'{a}bor Luk\'{a}cs for numerous comments and helpful suggestions.  
The author's collaboration on this
manuscript has started during the 49th Workshop ``Advances in Set-Theoretic Topology: Conference in Honour of Tsugunori Nogura on
his 60th Birthday'' of the International School of Mathematics ``G.~Stampacchia'' held on June 9--19, 2008 at the Center for Scientific
Culture ``Ettore Majorana'' in Erice, Sicily (Italy). The authors would like to express their warmest gratitude to the Ettore Majorana
Foundation and Center for Scientific Culture for providing excellent conditions which inspired this research endeavor.

\def\af{\sc}
\def\tf{\rm}
\def\papertitle#1{{#1}}
\def\collf{\em}
\def\bookf{\em}
\def\jf{\sl}

\end{document}